\documentclass[10pt]{article}
\usepackage[utf8]{inputenc}
\usepackage{tikz,amsfonts,amsmath}
\usepackage{natbib}
\usepackage{hyperref,doi}

\newcommand{\Gr}{\operatorname{Gr}}
\usepackage{xcolor,colortbl,tikz}
\usepackage[margin=1in]{geometry}

\title{Cluster Algebras and Symmetries of Regular Tilings}
\author{Adam Scherlis}

\begin{document}

\maketitle

\abstract{The classification of Grassmannian cluster algebras resembles that of regular polygonal tilings. We conjecture that this resemblance may indicate a deeper connection between these seemingly unrelated structures.}

\section{Introduction}

\definecolor{lightgreen}{rgb}{0.5,1,0.5}
\definecolor{lightred}{rgb}{1,0.5,0.5}
\newcommand{\redcell}[1]{\cellcolor{lightred}#1}
\newcommand{\greencell}[1]{\cellcolor{lightgreen}#1}
\newcommand{\yellowcell}[1]{\cellcolor{yellow}#1}
\newcolumntype{L}[1]{>{\raggedright\let\newline\\\arraybackslash\hspace{0pt}}m{#1}}
\newcolumntype{C}[1]{>{\centering\let\newline\\\arraybackslash\hspace{0pt}}m{#1}}
\newcolumntype{R}[1]{>{\raggedleft\let\newline\\\arraybackslash\hspace{0pt}}m{#1}}

Grassmannian cluster algebras and regular tilings (and the symmetry groups of those tilings) obey
similar classifications into finite, affine, and hyperbolic (or spherical, planar, and hyperbolic) cases. In fact, one can construct identical tables of cluster algebras and tilings, as shown below.\\

Formally, we observe that \textbf{the cluster algebra Gr($p,p+q$) is finite iff the Coxeter group $[p,q]$ is
finite, extended-affine iff $[p,q]$ is affine, and hyperbolic iff $[p,q]$ is
hyperbolic.}\\

Or, in simpler language: \textbf{Gr($p,p+q$) is of finite type iff the regular tiling $\{p,q\}$ is
spherical, infinite but of finite mutation type iff $\{p,q\}$ is planar, and of infinite mutation type iff $\{p,q\}$ is
hyperbolic.}\\

We conjecture that this similarity follows from a deeper connection between these structures. Understanding such a connection could provide a shorter, more elegant classification theorem for Grassmannian cluster algebras.

\section{Grassmannian Cluster Algebras}
\cite{math/0311148v1} and \cite{math/0608367v3}
have proven that the following classification of Grassmannian cluster algebras Gr($p,p+q$)
(with subspace dimension $p$ and codimension $q$) holds:
\begingroup
\renewcommand{\arraystretch}{2}
\begin{table}[h!]
\centering
\begin{tabular}{c|C{1.4cm}|C{1.4cm}|C{1.4cm}|C{1.4cm}|C{1.4cm}|C{1.4cm}|}
$_p \backslash ^q$ &2&3&4&5&6&7\\\hline
2&\greencell{$A_1$} & \greencell{$A_2$} & \greencell{$A_3$} & \greencell{$A_4$} & \greencell{$A_5$} & \greencell{$A_6$} \\\hline
3&\greencell{$A_2$} & \greencell{$D_4$} & \greencell{$E_6$} & \greencell{$E_8$} & \yellowcell{$E_8^{(1,1)}$} & \redcell{Gr(3,10)} \\\hline
4&\greencell{$A_3$} & \greencell{$E_6$} & \yellowcell{$E_7^{(1,1)}$} & \redcell{Gr(4,9)} & \redcell{Gr(4,10)} & \redcell{Gr(4,11)} \\\hline
5&\greencell{$A_4$} & \greencell{$E_8$} & \redcell{Gr(5,9)} & \redcell{Gr(5,10)} & \redcell{Gr(5,11)} & \redcell{Gr(5,12)} \\\hline
6&\greencell{$A_5$} & \yellowcell{$E_8^{(1,1)}$} & \redcell{Gr(6,10)} & \redcell{Gr(6,11)} & \redcell{Gr(6,12)} & \redcell{Gr(6,13)} \\\hline
7&\greencell{$A_6$} & \redcell{Gr(7,10)} & \redcell{Gr(7,11)} & \redcell{Gr(7,12)} & \redcell{Gr(7,13)} & \redcell{Gr(7,14)} \\\hline
\end{tabular}
\end{table}\endgroup
\begin{itemize}\item Green cells of the table denote finite cluster algebras, with finite Dynkin diagrams $X_n$.
\item Yellow cells are infinite, but of finite mutation type, with extended affine Dynkin diagrams $X_n^{(1,1)}$.
\item Red cells are of infinite mutation type, with hyperbolic Dynkin diagrams.\end{itemize}

This classification depends on the parameter $r=(p-2)(q-2)$ (see \cite{math/0608367v3} Prop. 12.11): \textbf{Gr($p,p+q$) is finite for $r<4$, infinite with finite mutation type for $r=4$, and infinite mutation type for $r>4$.}\\

\cite{math/0608367v3} prove this statement by considering individual cases.\\

There is a canonical isomorphism Gr($p,p+q$)$\cong$Gr($q,q+p$).

\section{Regular Tilings}
This table depicts regular tilings $\{p,q\}$. These are two-dimensional surfaces formed by joining together regular $p$-gons, with $q$ at each vertex.
\begingroup
\renewcommand{\arraystretch}{2}
\begin{table}[h!]
\centering
\begin{tabular}{c|C{1.4cm}|C{2.2cm}|C{2.2cm}|C{2.2cm}|C{2.2cm}|C{1.4cm}|}
$_p \backslash ^q$ &2&3&4&5&6&7\\\hline
2&\greencell{\{2,2\}} & \greencell{\{2,3\}} & \greencell{\{2,4\}} & \greencell{\{2,5\}} & \greencell{\{2,6\}} & \greencell{\{2,7\}} \\\hline
3&\greencell{\{3,2\}} & \greencell{tetrahedron} & \greencell{octahedron} & \greencell{icosahedron} & \yellowcell{triangular tiling} & \redcell{\{3,7\}} \\\hline
4&\greencell{\{4,2\}} & \greencell{cube} & \yellowcell{square tiling} & \redcell{\{4,5\}} & \redcell{\{4,6\}} & \redcell{\{4,7\}} \\\hline
5&\greencell{\{5,2\}} & \greencell{dodecahedron} & \redcell{\{5,4\}} & \redcell{\{5,5\}} & \redcell{\{5,6\}} & \redcell{\{5,7\}} \\\hline
6&\greencell{\{6,2\}} & \yellowcell{hexagonal tiling} & \redcell{\{6,4\}} & \redcell{\{6,5\}} & \redcell{\{6,6\}} & \redcell{\{6,7\}} \\\hline
7&\greencell{\{7,2\}} & \redcell{\{7,3\}} & \redcell{\{7,4\}} & \redcell{\{7,5\}} & \redcell{\{7,6\}} & \redcell{\{7,7\}} \\\hline
\end{tabular}
\end{table}\endgroup\\
\begin{itemize}\item Green cells are spherical tilings: hosohedra $\{2,q\}$, dihedra $\{p,2\}$, and the five Platonic solids.
\item The three yellow cells are planar tilings.
\item The red cells are regular hyperbolic tilings.\end{itemize}

The nature of such a tiling depends on the value of
$r=(p-2)(q-2)$: it will be spherical for $r<4$, planar for $r=4$, and hyperbolic for $r>4$. This can be shown by calculating the angular defect at each vertex.\\

The tilings $\{p,q\}$ and $\{q,p\}$ are dual.

\subsection{Symmetry Groups of Tilings}
The tiling $\{p,q\}$ has the symmetry group $[p,q]$ in Coxeter notation; these have associated Coxeter-Dynkin diagrams and obey a Cartan-Killing classification, much like cluster algebras.
\begin{itemize}
\item Spherical tilings correspond to finite Coxeter groups $X_n$.
\item Planar tilings have affine Coxeter groups $X_n^{(1)}$ (also called $\widetilde X_n$).
\item Hyperbolic tilings have hyperbolic Coxeter groups.\end{itemize}
Dual tilings have the same symmetries, $[p,q]\cong [q,p]$.

\section{Dynkin diagrams}
Note that $[p,q]$ and $\Gr(p,p+q)$ have very different, and apparently unrelated, Dynkin diagrams. In particular, the diagram for $[p,q]$ is of rank 3 while $\Gr(p,p+q)$ has rank $(p-1)(q-1)$. The initial quiver of $\Gr(p,p+q)$ is mutation-equivalent to its Dynkin diagram, which is also of this rank.
\begin{figure}[h!]
\centering
\begin{minipage}{.5\textwidth}
  \centering
  \includegraphics[width=.4\linewidth]{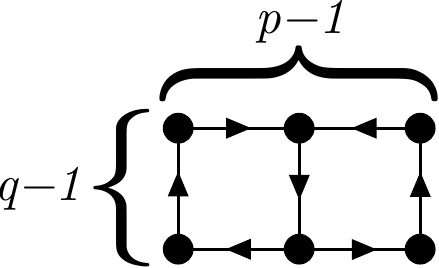}
  \caption{The initial quiver of $\Gr(p,p+q)$}
  \label{fig:test1}
\end{minipage}%
\begin{minipage}{.5\textwidth}
  \centering
  \includegraphics[width=.4\linewidth]{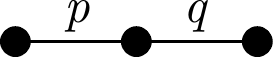}
  \caption{The Coxeter-Dynkin diagram of $[p,q]$}
  \label{fig:test2}
\end{minipage}
\end{figure}

\section{Summary Table}
\begin{tabular}{r||l|l||l|l|l}
  $r$ & Cluster Algebra &  & Regular Tiling& Symmetries  & Coxeter Notation\\\hline\hline
$r=0$ &\textbf{simplest finite type} & &  \textbf{degenerate spherical}&\\
&$\Gr(2,p+2)\cong \Gr(p,p+2)$ & $A_{p-1}$ & hosohedron $\cong$ dihedron& $A_1\times I_2(p)$  & $[2,p]\cong [p,2]$\\
\hline
$0<r<4$ &\textbf{other finite type} &  & \textbf{Platonic}&\\
& $\Gr(3,6)$ & $D_4$ &  tetrahedron & $A_3$ &$[3,3]$\\
 & $\Gr(3,7)\cong\Gr(4,7)$ & $E_6$ & octahedron $\cong$ cube& $BC_3$  & $[3,4]\cong [4,3]$\\
 & $\Gr(3,8)\cong\Gr(5,8)$ & $E_8$ &  icosahedron $\cong$ dodecahedron &$H_3$ & $[3,5]\cong [5,3]$\\
\hline
$r=4$ &\textbf{finite mutation type} &  & \textbf{planar}&\\
& $\Gr(4,8)$ & $E_7^{(1,1)}$ &  square tiling &$C_2^{(1)}$ & $[4,4]$\\
& $\Gr(3,9)\cong\Gr(6,9)$ & $E_8^{(1,1)}$ & triangular $\cong$ hexagonal tiling &$G_2^{(1)}$ &  $[3,6]\cong[6,3]$\\[0.2em]
\hline
$r>4$ &\textbf{infinite mutation type} &  & \textbf{hyperbolic}&
\end{tabular}

\section{Acknowledgements}
The author would like to thank D. Parker, M. Spradlin, A. Volovich, A. Goncharov, and M. Shapiro for enlightening discussions.

\bibliographystyle{plainnat}
\bibliography{references}

\end{document}